\documentclass[11pt]{article}%
\usepackage{amsmath}
\usepackage{amsfonts}
\usepackage{amssymb}
\usepackage{graphicx}
\usepackage{a4wide}
\usepackage{xcolor}%
\providecommand{\U}[1]{\protect\rule{.1in}{.1in}}
\newtheorem{theorem}{Theorem}

\newtheorem{remark}[theorem]{Remark}

\begin{document}

\title{On the definition of \textquotedblleft almost LUR (ALUR)" notion}
\author{Constantin Z\u{a}linescu\thanks{Octav Mayer Institute of Mathematics, Ia\c{s}i
Branch of Romanian Academy, Ia\c{s}i, Romania, email:
\texttt{zalinesc@uaic.ro}}}
\date{}
\maketitle

\begin{abstract}
The notion of almost LUR (ALUR) point is introduced in the paper
[P.\ Bandyopadhyay et al., Some generalizations of locally uniform rotundity,
J.~Math.\ Anal.\ Appl., 252, 906--916 (2000)], where one says that the point
$x$ of the unit sphere $S_{X}$ of a Banach space is an almost LUR (ALUR) point
of $B_{X}$ if for any sequences $\{x_{n}\}\subseteq B_{X}$ and $\{x_{m}^{\ast
}\}\subseteq B_{X^{\ast}}$, the condition $\lim_{m}\lim_{n}x_{m}^{\ast}\left(
\frac{x_{n}+x}{2}\right)  =1$ implies $\lim_{n}\left\Vert x_{n}-x\right\Vert
=0$, without mentioning what is meant by $\lim_{m}\lim_{n}\gamma_{m,n}=\gamma$
for $\gamma$, $\gamma_{m,n}\in\mathbb{R}$; $X$ is ALUR if $X$ is almost LUR at
any $x\in$ $S_{X}$. Of course, the natural definition for this iterated limit
would be that for each $m$ sufficiently large there exists $\mu_{m}%
:=\lim_{n\rightarrow\infty}\gamma_{m,n}\in\mathbb{R}$ and $\gamma
=\lim_{m\rightarrow\infty}\mu_{m}$. However, as seen in some proofs
where $\lim_{m}\lim_{n}$ appears, this interpretation is not
confirmed. In this paper we examine several works in which almost
LUR is mentioned and, especially, the proofs of those results in
which the above definition of \textquotedblleft almost LUR" point
(or space) is invoked. Moreover, we analyze similar problems related
to the notion CWALUR which extend ALUR. Furthermore, we mention
several gaps in the proofs of some results. Finally, we propose the
change of $\lim_{m}\lim_{n}$ by $\lim_{m}\liminf_{n}$ in the
definitions of several types of ALUR points; moreover, we provide
the complete proofs of two results from the literature in which
$\lim_{m}\lim_{n}$ were used effectively, using
$\lim_{m}\liminf_{n}$ instead.

\end{abstract}

\bigskip

P.\ Bandyopadhyay, D. Huang, B.-L.\ Lin and S.~L.\ Troyanski say, in
\cite[Def.\ 2]{BaHuLiTr00},\footnote{In \cite{Hua99} one finds the following
text: \textquotedblleft Definition 5.1.1.\ [16] Let $X$ be a Banach space,
$x\in S_{X}$ is called an almost locally uniformly rotund (ALUR in short)
(resp.\ weakly almost locally uniformly rotund (WALUR in short)) point of
$B_{X}$ if for any $x_{n}\in B_{X}$, $n\in\mathbb{N}$ and $x_{m}^{\ast}\in
S_{X^{\ast}}$, $m\in\mathbb{N}$ such that $\lim_{m}\lim_{n}x_{m}^{\ast}%
(x_{n})=\lim_{m}x_{m}^{\ast}(x)=1$, then $\lim_{n}x_{n}=x$ (resp.\ $w$%
-$\lim_{n}x_{n}=x$). $X$ is said to be ALUR (resp.\ WALUR) if every $x\in
S_{X}$ is an ALUR (resp.\ WALUR) point of $B_{X}$." ([16] seems to contain a
previously submitted version of \cite{BaHuLiTr00}.)} that the
\textquotedblleft point $x$ of the unit sphere $S_{X}$ of a Banach space $X$
... is an an almost LUR (ALUR) (resp.\ weakly almost LUR (WALUR)) point of
$B_{X}$ if for any $\{x_{n}\}\subseteq B_{X}$ and $\{x_{m}^{\ast}\}\subseteq
B_{X^{\ast}}$, the condition $\lim_{m}\lim_{n}x_{m}^{\ast}\left(  \frac
{x_{n}+x}{2}\right)  =1$ implies $\lim_{n}\left\Vert x_{n}-x\right\Vert =0$
(resp.\ w-$\lim_{n}(x_{n}-x)=0$)"; in \cite[Def.\ 5]{BaHuLiTr00} one says
\textquotedblleft that $x^{\ast}\in S_{X^{\ast}}$ is a w*-ALUR point of
$B_{X^{\ast}}$ if for any $\{x_{n}^{\ast}\}\subseteq B_{X^{\ast}}$ and
$\{x_{m}\}\subseteq B_{X}$ the condition $\lim_{m}\lim_{n}\left(  \frac
{x_{n}^{\ast}+x^{\ast}}{2}\right)  (x_{m})=1$ implies w*-$\lim x_{n}^{\ast
}=x^{\ast}$".

\smallskip Unfortunately, it is not defined (explained) what is meant by
$\lim_{m}\lim_{n}\gamma_{m,n}=\gamma$ when $\gamma,$ $\gamma_{m,n}%
\in\mathbb{R}$ (for $m,n\in\mathbb{N}^{\ast}:=\{1,2,...\}$). Our guess was
that the preceding equality is the one formulated in the next definition.

\smallskip\textbf{Definition A}. Let $\gamma,$ $\gamma_{m,n}\in\mathbb{R}$ for
$m,n\in\mathbb{N}^{\ast}$. One says that $\gamma=\lim_{m}\lim_{n}\gamma_{m,n}$
if for every $m\in\mathbb{N}^{\ast}$ there exists $\mu_{m}:=\lim_{n}%
\gamma_{m,n}$ and $\lim_{m}\mu_{m}=\gamma$.\footnote{On 24.11.2024 we asked
the following question to the authors of the paper \cite{BaHuLiTr00}:
\textquotedblleft Taking $\gamma_{n,m}:=x_{m}^{\ast}\left(  \frac{x_{n}+x}%
{2}\right)  ,$ does $\lim_{m}\lim_{n}\gamma_{n,m}=1$ mean that for each
$m\geq1$ there exists $\mu_{m}:=\lim_{n}\gamma_{n,m}$ and $\lim_{m}\mu_{m}%
=1?$"
After two days we received the following answer from one of the authors of
\cite{BaHuLiTr00}: \textquotedblleft I think your interpretation is what we
had in mind".}

\smallskip This interpretation seems to be confirmed by the proof of the
implication (a)~$\Rightarrow$ (b) of \cite[Th.\ 6]{BaHuLiTr00}, where
\textquotedblleft(a) $x^{\ast}$ is a rotund point of $B_{X^{\ast}}$; (b)
$x^{\ast}$ is a w*-ALUR point of $B_{X^{\ast}}$". Recall that, by
\cite[Def.\ 2(a)]{BaHuLiTr00}, $x\in S_{X}$ is \textquotedblleft a rotund
point of the unit ball $B_{X}$ of $X$ (or, $X$ is rotund at $x$) if
$\left\Vert y\right\Vert =\left\Vert (x+y)/2\right\Vert =1$ implies $x=y$".
Below we quote the announced proof:

\smallskip

\textquotedblleft(a)~$\Rightarrow$ (b). Let $\{x_{n}^{\ast}\}\subseteq
B_{X^{\ast}}$ and $\{x_{m}\}\subseteq B_{X}$ such that
\begin{equation}
\lim_{m}\lim_{n}\left(  \frac{x_{n}^{\ast}+x^{\ast}}{2}\right)  (x_{m})=1.
\label{rz1}%
\end{equation}
As $\{x_{n}^{\ast}\}\subseteq B_{X^{\ast}},$ $\{x_{n}^{\ast}\}$ has w*-cluster
points in $B_{X^{\ast}}$. Let $y^{\ast}$ be a w*-cluster point of
$\{x_{n}^{\ast}\}$. Then $\left\Vert y^{\ast}\right\Vert \leq1$. It follows
that
\[
\lim_{m}\left(  \frac{x^{\ast}+y^{\ast}}{2}\right)  (x_{m})=1
\]
and hence $\left\Vert (x^{\ast}+y^{\ast})/2\right\Vert =1$. By (a), we
conclude that $y^{\ast}=x^{\ast}$. This implies the sequence $\{x_{n}^{\ast
}\}$ has a unique w*-cluster point $x^{\ast}$, that is, w*-$\lim x_{n}^{\ast
}=x^{\ast}$." (The number is introduced by us.)

\smallskip To the above proof, we add some details. By (\ref{rz1}) and
Definition A, for $m\geq1$ there exists $\mu_{m}:=\lim_{n}\left(  \frac
{x_{n}^{\ast}+x^{\ast}}{2}\right)  (x_{m})$ and $\lim_{m}\mu_{m}=1$. Because
$y^{\ast}$ is a w*-cluster point of $\{x_{n}^{\ast}\}$, {there exists a subnet
$\{x_{n_{i}}^{\ast}\}$ {of} $\{x_{n}^{\ast}\}$ {such that }$x_{n_{i}}^{\ast
}\rightarrow^{w^{\ast}}y^{\ast}$ $(\in B_{X^{\ast}})$, whence $\left\Vert
y^{\ast}\right\Vert \leq1$ and so
\[
\mu_{m}=\lim_{i}\left(  \frac{x_{n_{i}}^{\ast}+x^{\ast}}{2}\right)
(x_{m})=\left(  \frac{y^{\ast}+x^{\ast}}{2}\right)  (x_{m})\leq\tfrac{1}%
{2}\left\Vert x_{m}\right\Vert \cdot\left\Vert y^{\ast}+x^{\ast}\right\Vert
\leq\tfrac{1}{2}\left\Vert y^{\ast}+x^{\ast}\right\Vert \leq1
\]
for $m\geq1$. Because $\lim_{m}\mu_{m}=1$ one obtains that $\tfrac{1}%
{2}\left\Vert y^{\ast}+x^{\ast}\right\Vert =1$. Using (a) one gets $y^{\ast
}=x^{\ast}$. This shows that $x^{\ast}$ is the unique w*-cluster point of the
sequence }$\{x_{n}^{\ast}\}${, and so w*-$\lim x_{n}^{\ast}=x^{\ast}$. Hence
$x^{\ast}$ is a w*-ALUR point of }$B_{X^{\ast}}$.

\medskip

Let us have a look to the proof of the implication (b)~$\Rightarrow$ (c) of
\cite[Th.\ 6]{BaHuLiTr00}, where (b) was recalled above and \textquotedblleft%
(c) for every unbounded nested sequence $\{B_{n}\}$ of balls such that
$x^{\ast}$ is bounded below on $\cup B_{n}$, if for any $\{y_{n}^{\ast
}\}\subseteq S_{X^{\ast}}$, the sequence $\{\inf y_{n}^{\ast}(B_{n})\}$ is
bounded below, then w*-$\lim y_{n}^{\ast}=x^{\ast}$".\footnote{In the sequel
we recall only those notions and notations which are used in our arguments.}
The proof of the implication (b)~$\Rightarrow$ (c) ends with the following text:

\smallskip\textquotedblleft If we now put $y_{n}=x_{n}/r_{n}$, it follows that
$\left\Vert y_{n}\right\Vert <1$. The fact that $\{B_{n}\}$ is nested implies
that $y_{n}^{\ast}(y_{m})\geq1+c/r_{m}$ for all $n\geq m$. It follows that
\[
\left(  \frac{x^{\ast}+y_{n}^{\ast}}{2}\right)  (y_{m})\geq1+c/r_{m}%
\]
for all $n\geq m$. Since $r_{m}\rightarrow\infty$, we conclude%
\[
\lim_{m}\lim_{n}\left(  \frac{x^{\ast}+y_{n}^{\ast}}{2}\right)  (y_{m})=1.
\]
Hence, w*-$\lim y_{n}^{\ast}=x^{\ast}$."

\smallskip

Clearly, only the fact that $\left(  \frac{x^{\ast}+y_{n}^{\ast}}{2}\right)
(y_{m})\geq1+c/r_{m}$ \emph{for }$n\geq m\geq1$\emph{ does not imply that
}$\lim_{n}\left(  \frac{x^{\ast}+y_{n}^{\ast}}{2}\right)  (y_{m})$\emph{
exists} in $\mathbb{R}$.

\medskip

So, \emph{which is the meaning of the equality }$\lim_{m}\lim_{n}x_{m}^{\ast
}\left(  \frac{x_{n}+x}{2}\right)  =1$ from \cite[Def.\ 2]{BaHuLiTr00} or,
more generally, \emph{the meaning of} $\lim_{m}\lim_{n}\gamma_{m,n}=\gamma$?

\medskip We found other 16 works in which the same definition of ALUR (with
the sense ``almost LUR'') point (or space) is given; see \cite[Def.\ 2.3(ii)]%
{BanLin01}, \cite[Def.\ 1.2(iii)]{BaHuLi04}, \cite[p.\ 132]{AizGar05},
\cite[Def.\ 2.19]{BanGod06}, \cite[Def.\ 3.1]{BaLiLiNa08}, \cite[Def.\ 2.5(b)]%
{BaLiRa09}, \cite[Def.\ 7.2.11(c)]{Pau10}, \cite[Def.\ 6.6(c)]{BanPau10},
\cite[p.\ 3897]{ZhaLiu11}, \cite[p.\ 2]{ZhaLiu11b}, \cite[p.\ 2]{ZhaLiu12},
\cite[p.\ 285]{DamBaj17}, \cite[p.\ 114]{GupNar17}, \cite[Def.\ 1.1(2)]%
{ZhZhLi17}, \cite[p.\ 20]{JuMaRu23}, \cite[p.\ 1692]{DeBSom24}.

Moreover, in \cite[Def.\ 1.10(b)]{DasPau24}, one considers $\lim_{m,n}%
h_{m}\left(  \frac{x_{n}+x}{2}\right)  =1$ with $h_{m}\in S_{X^{\ast}}$ when
defining \textquotedblleft$\tau$-almost locally uniformly rotund ($\tau$-ALUR
in short)" Banach spaces. Is $\lim_{m,n}h_{m}\left(  \frac{x_{n}+x}{2}\right)
=1$ the same as $\lim_{m}\lim_{n}h_{m}\left(  \frac{x_{n}+x}{2}\right)  =1$?

Notice that even if the string \textquotedblleft almost LUR" does appear four
times in \cite[p.\ 2]{KaLoMaWe20}, the corresponding definition is not
provided, \cite[Prop.\ 5.3]{KaLoMaWe20} being attributed to \cite[Prop.\ 2.23]%
{BanGod06}.

\medskip

We hoped to deduce the meaning of $\lim_{m}\lim_{n}\gamma_{m,n}=\gamma$ by
following, in the above works, the proofs in which ALUR (having the meaning of
almost LUR) is present. Of course, the proofs of those results have to contain
an iterated limit ($\lim_{m}\lim_{n}$), and so at least one limit ($\lim$).

\smallskip We already mentioned the implication (b)~$\Rightarrow$ (c) of
\cite[Th.\ 6]{BaHuLiTr00} in which ALUR appears as hypothesis and the
existence of $\lim_{n}\left(  \frac{x^{\ast}+y_{n}^{\ast}}{2}\right)  (y_{m})$
is not proved (for $m\geq1$). Other statements from \cite{BaHuLiTr00} in which
ALUR appears are: the diagram on page 908 for which one considers that
\textquotedblleft It is necessary to prove only ALUR~$\rightarrow$ MLUR and
WALUR~$\rightarrow$ WMLUR",

\smallskip\textquotedblleft Corollary 8. Let $X$ be a Banach space. For $x\in
S_{X}$, the following are equivalent:

(a) $x$ is a rotund point of $B_{X^{\ast\ast}}$;

(b) $x$ is a w*-ALUR point of $B_{X^{\ast\ast}}$;

(b$^{\prime}$) $x$ is a wALUR point of $B_{X}$;

(c) for every unbounded nested sequence $\{B_{n}^{\ast}\}$ of balls in
$X^{\ast}$ such that $x$ is bounded below on $\cup B_{n}^{\ast}$, if for any
$\{y_{n}^{\ast\ast}\}\subseteq S_{X^{\ast\ast}}$, the sequence $\{\inf
y_{n}^{\ast\ast}(\cup B_{n}^{\ast})\}$ is bounded below, then w*-$\lim
y_{n}^{\ast\ast}=x$;

(c$^{\prime}$) ...; (d) ...; (e) ... on $\cup B_{n}^{\ast}$, $\cup B_{n}%
^{\ast}$ is an affine half-space determined by $x$.",

\smallskip\noindent and \textquotedblleft Corollary 12. There exists an ALUR
Banach space which fails to be WLUR.".

\smallskip In what concerns \cite[Cor.\ 8]{BaHuLiTr00} one finds

\textquotedblleft Proof. The equivalence of (a) to (e) is Theorem 6, while
(b)~$\Rightarrow$ (b$^{\prime}$) and (c)~$\Rightarrow$ (c$^{\prime}$) are
immediate. And (b$^{\prime}$)~$\Rightarrow$ (c$^{\prime}$) follows similarly
as (b)~$\Rightarrow$ (c) of Theorem 6. (c$^{\prime}$)~$\Rightarrow$ (d). Let
$\{B_{n}^{\ast}=B(x_{n}^{\ast},r_{n})\}$ ... a contradiction."

\smallskip So, in the proof of \cite[Cor.\ 8]{BaHuLiTr00}, (b) and
(b$^{\prime}$) are not involved directly; in fact, $\lim$ does not appear at
all in that proof.

\smallskip Concerning \cite[Cor.\ 12]{BaHuLiTr00} one finds:

\smallskip\textquotedblleft Proof. Let $(X,\left\vert \cdot\right\vert )$ be
the space from Proposition 11 which fails to be WLUR. By Corollary 8,
$(X,\left\vert \cdot\right\vert )$ is WALUR. Since $(X,\left\vert
\cdot\right\vert )$ has the KK property, we get that $(X,\left\vert
\cdot\right\vert )$ is ALUR. $\square$\textquotedblright

\smallskip Hence, in the proofs of Corollaries 8 and 12 the ALUR definitions
are not used.

\smallskip Let us analyze the proof of the implications \textquotedblleft
ALUR~$\rightarrow$ MLUR and WALUR~$\rightarrow$ WMLUR". First, by
\cite[Def.\ 2 (d)]{BaHuLiTr00}, $x\in S_{X}$ is \textquotedblleft a midpoint
LUR (MLUR) (resp.\ weakly midpoint LUR (WMLUR)) point of $B_{X}$ if $\lim
_{n}\left\Vert x\pm x_{n}\right\Vert =1$ implies $\lim_{n}\left\Vert
x_{n}\right\Vert =0$ (resp.\ w-$\lim_{n}x_{n}=0$)".

The provided proof is the following:

\smallskip\textquotedblleft Indeed let $x\in S_{X}$ be an ALUR (resp.\ WALUR)
point of $B_{X}$ and $\lim_{n}\left\Vert x\pm x_{n}\right\Vert =1$. Then for
every $x^{\ast}\in S_{X^{\ast}}$ such that $x^{\ast}(x)=1$, we have $\lim
_{n}x^{\ast}(x_{n})=0$. So $\lim_{n}x^{\ast}\left(  (x+(x_{n}+x))/2\right)
=1$, hence $\lim_{n}\left\Vert x_{n}\right\Vert =0$ (resp.\ w-$\lim_{n}%
x_{n}=0$)."

\smallskip We agree that taking $\{x_{n}\}\subseteq X$ such that $\lim
_{n}\left\Vert x\pm x_{n}\right\Vert =1$ and $x^{\ast}\in S_{X^{\ast}}$ such
that $x^{\ast}(x)=1$ one gets $\lim_{n}x^{\ast}(x_{n})=0$, and so $\lim
_{n}x^{\ast}\left(  (x+(x_{n}+x))/2\right)  =1$.

\medskip

One must understand that one applies the fact that $X$ is ALUR (resp.\ WALUR)
at $x$ $(\in S_{X})$. So, one has to identify some sequences $\{x_{m}^{\ast
}\}\subseteq B_{X^{\ast}}$ and $(x_{n}^{\prime})\subseteq B_{X}$ verifying
{$\lim_{m}\lim_{n}x_{m}^{\ast}\left(  \frac{x_{n}^{\prime}+x}{2}\right)  =1$
(in order to get }$x_{n}^{\prime}\rightarrow x$). Observe that $\lim
_{n}x^{\ast}\left(  (x+(x_{n}+x))/2\right)  =1$ can be written as
\[
{\lim_{n}x^{\ast}\left(  \frac{(x_{n}+x)+x}{2}\right)  }=\lim_{n}x^{\ast
}\left(  \frac{(x+(x_{n}+x))}{2}\right)  =\lim_{n}x^{\ast}\left(
(x+(x_{n}+x))/2\right)  {=1}.
\]
So, naturally, one takes $x_{m}^{\ast}:=x^{\ast}$ $(\in S_{X^{\ast}}\subseteq
B_{X^{\ast}})$ for $m\geq1$ and {$x_{n}^{\prime}:=x_{n}+x$} for $n\geq1$, and
so {$\lim_{m}\lim_{n}x_{m}^{\ast}\left(  \frac{x_{n}^{\prime}+x}{2}\right)
=1$ in the sense of Definition A. Is $\{x_{n}^{\prime}\}\subseteq B_{X}$ as in
Definition 2(c)?} What we know is that $\lim_{n}\left\Vert x\pm x_{n}%
\right\Vert =1$ (probably with $\{x_{n}\}\subseteq X$ (or, maybe,
$\{x_{n}\}\subseteq B_{X}$ as in Definition 2(b)\&(c))). Possibly, $\{x_{n}\}$
could be an arbitrary sequence from $X$ in the definition of ALUR (but this
must be proved or, at least, mentioned)!

\medskip

Notice that $\lim_{m}\lim_{n}$ appears 4 times in \cite{BaHuLiTr00} (including
Definitions 1 and 5) and $\lim$ 46 times.

\medskip Let us continue to look in the works mentioned above for the proofs
of the implications $p\Rightarrow q$ in which ALUR (with the sense ``almost
LUR'') is involved in the hypothesis or the conclusion.

\medskip

-- \cite{BanLin01}: Notice that $\lim$ appears 22 times including $\lim
_{m}\lim_{n}$ 2 times (in Definition 2.3(ii)). Observe that \cite[Th.\ 6 and
Cor.\ 8]{BaHuLiTr00} are quoted as \cite[Ths.\ 2.4 and 2.6]{BanLin01},
respectively. Moreover, only Theorem 2.1 (which asserts that conditions (a)
and (e) of Theorem 2.4 are equivalent) and Theorem 2.5 have proofs, but they
do not refer to ALUR. All the other results and remarks are included in
\cite{BaHuLi04}, where the majority of the results have proofs.

\medskip

-- \cite{BaHuLi04}: Notice that $\lim$ appears 79 times, including $\lim
_{m}\lim_{n}$ 3 times (2 in Definition 2.3(ii) and 1 in the proof of
Proposition 2.5).

First observe that \cite[Th.\ 6 and Cor.\ 8]{BaHuLiTr00} are quoted as
\cite[Th.\ 1.4 and 1.5]{BaHuLi04}, respectively.

The proof of the next result use the definition of an ALUR point in a similar
manner to \cite[Th.\ 6 (b)~$\Rightarrow$ (c)]{BaHuLiTr00}.

\smallskip\textquotedblleft{Proposition} 2.5 For a Banach space $X$ and
$x^{\ast}\in S(X^{\ast})$, the following are equivalent:

(a) $x^{\ast}$ is a rotund point of $B(X^{\ast})$

(b) For any $\{x_{n}^{\ast}\}\subseteq B(X^{\ast})$, if $\{(x_{n}^{\ast
}+x^{\ast})/2\}$ is asymptotically normed by $B(X)$, then w*-$\lim x_{n}%
^{\ast}=x^{\ast}$.

\smallskip

{Proof} (a) $\Rightarrow$ (b). Suppose $\{y_{n}^{\ast}\}\subseteq B(X^{\ast})$
is asymptotically normed by $B(X)$. By definition, {for any $m\geq1$, there
exists a $x_{m}\in B$}$(X)${ and $N_{m}\in\mathbb{N}$ such that $y_{n}^{\ast
}(x_{m})>1-1/m$ for all $n\geq N_{m}$. It follows that $\lim_{m}\lim_{n}%
y_{n}^{\ast}(x_{m})=1$.} Thus, {(a) $\Rightarrow$ (b) follows from Theorem 1.4
(b)}.

(b) $\Rightarrow$ (a). We prove (b) $\Rightarrow$ Theorem 1.4 (c). ...''

\smallskip The proof does not contain other $\lim$. Even if not said
explicitly, {$y_{n}^{\ast}$ is $(x_{n}^{\ast}+x^{\ast})/2$ in the above proof.
By the implication (a)~}$\Rightarrow$ (b) of Theorem 1.4 (that is,
\cite[Th.\ 6]{BaHuLiTr00}), $X^{\ast}$ is w*-ALUR at $x^{\ast}$. Then, because
{$\lim_{m}\lim_{n}\left(  \frac{x_{n}^{\ast}+x^{\ast}}{2}\right)  (x_{m}%
)=\lim_{m}\lim_{n}y_{n}^{\ast}(x_{m})=1$ one gets }w*-$\lim x_{n}^{\ast
}=x^{\ast}${ by }Definition 1.2(iii).

\smallskip But \emph{which is the definition of }$\lim_{m}\lim_{n}\gamma
_{m,n}=\gamma$\emph{ for getting }$\lim_{m}\lim_{n}y_{n}^{\ast}(x_{m}%
)=1$\emph{ from the fact that }$y_{n}^{\ast}(x_{m})>1-1/m$\emph{ for all
}$n\geq N_{m}$\emph{ (for all $m\geq1$\emph{)}}?

\smallskip Having in view the above discussion and the fact that Theorem 1.4
(b) is only a logical proposition (which might be true or false), possibly,
one wished to say {\textquotedblleft Theorem 1.4 (a)$~\Rightarrow$ (b)"
instead of \textquotedblleft Theorem 1.4 (b)".}

\smallskip In the sequel we mention those results from the works in our list
which contain assertions (conditions) referring to versions of ALUR (as
abbreviation for almost LUR) in their formulation and/or their proofs, or
iterated limits, as well as the number of times $\lim$ and $\lim_{m}\lim_{n}$
do appear in the proofs.

\smallskip We mention below the other results from \cite{BaHuLi04} referring
to ALUR points and the corresponding proofs concerning them.

\smallskip\textquotedblleft{Proposition} 4.4 For a Banach space $X$ and $x\in
S(X)$, the following are equivalent:

(a) $x$ is a wALUR point of $B(X)$. (b) ...

(c) $x$ is w*-exposed in $B(X^{\ast\ast})$ by every $x^{\ast}\in D(x)$. (d)
... (e) ... (f) ...

(g) For every $x^{\ast}\in D(x)$ and for any $\{x_{n}\}\subseteq S(X)$, if
$x^{\ast}(x_{n})\rightarrow1$, then w-$\lim x_{n}=x$.",

\noindent where, \textquotedblleft$D(x)=\{x^{\ast}\in S(X^{\ast}):x^{\ast
}(x)=1\}$, $x\in S(X)$" (cf.\ \cite[Def.\ 1.13(a)]{BaHuLi04}).

\smallskip

The proof contains one $\lim$. The parts of the proof referring to (a) are the
following: \textquotedblleft(a)~$\Leftrightarrow$ (b) is similar to
Proposition 2.5. (a)~$\Rightarrow$ (c) follows from Lemma 4.1", where

\textquotedblleft{Lemma} 4.1 Let $X$ be a Banach space. $x\in S(X)$ is a
rotund point of $B(X)$ if and only if $x$ is exposed by every $x^{\ast}\in
D(x)$."

\textquotedblleft(g)~$\Rightarrow$ (a). We prove (g)~$\Rightarrow$ Theorem
1.5(c$^{\prime}$)". As mentioned above, Theorem 1.5 is \cite[Cor.\ 8]%
{BaHuLiTr00}).

\smallskip

\textquotedblleft{Corollary} 4.6 For a Banach space $X$ and $x\in S(X)$, the
following are equivalent:

(a) $x$ is a ALUR point of $B(X)$.

(b) For any $\{x_{n}\}\subseteq B(X)$, if $\{(x_{n}+x)/2\}$ is asymptotically
normed by $B(X^{\ast})$, then $\lim x_{n}=x$.

(c) $x$ is a wALUR point as well as a PC of $B(X)$.

(d) $x$ is a wALUR point as well as a seq PC of $B(X)$. (e) ... (f) ... (g) ...

(h) for every unbounded nested sequence $\{B_{n}^{\ast}\}$ of balls in
$X^{\ast}$ such that $x$ is bounded below on $\cup B_{n}^{\ast}$, if for any
$\{y_{n}\}\subseteq S(X)$, the sequence $\{\inf y_{n}(B_{n}^{\ast})\}$ is
bounded below, then $\lim y_{n}=x$."

\smallskip The proof does not contain $\lim$. The parts of the proof referring
to (a), (c) and (d) are the following: \textquotedblleft(a)\ $\Leftrightarrow$
(b) is similar to Theorem 2.5. (c)~$\Rightarrow$ (d)~$\Rightarrow$
(a)~$\Rightarrow$ (e) is immediate from the definitions"; \textquotedblleft%
(h)~$\Rightarrow$ (c). By Theorem 1.5, (h) implies $x$ is a wALUR point of
$B(X)$." Notice that there is no Theorem 2.5 in \cite{BaHuLi04}; possibly one
wished to refer to \cite[Prop.\ 2.5]{BaHuLi04} (see above) whose condition (b)
is similar to the present (b).

\smallskip\textquotedblleft{Proposition} 5.7 For a Banach space $X$ and $x\in
S(X)$, the following are equivalent:

I. (a) $x$ is a wALUR point of $B(X)$. (b) ...

II. (a) $x$ is an ALUR point of $B(X)$. (b) ..."

{Proof} This follows from Lemma 5.2, Proposition 4.4(g) and Corollary 4.6(e).
$\square$"

\smallskip Clearly, the proof does not contain $\lim$. Although, notice that
Proposition 4.4(g) and Corollary 4.6(e) are just two logical propositions
which might be true or false!

\smallskip

\textquotedblleft{Corollary} 5.9 Let $X$ be a Banach space.

I. (a) $X$ is wALUR if and only if every $x^{\ast}\in D(S(X))$ is a smooth
point of $B(X^{\ast})$. In particular, if $X^{\ast}$ is smooth, then $X$ is wALUR.

(b) $X$ is ALUR if and only if every $x^{\ast}\in D(S(X))$ is a Fr\'{e}chet
smooth point of $B(X^{\ast})$. In particular, if $X^{\ast}$ is Fr\'{e}chet
smooth, then $X$ is ALUR.

II. (a) ... (b) If wALUR points of $B(X^{\ast})$ form a boundary for $X$ (in
particular, if $X^{\ast}$ is wALUR), then $X$ is very smooth.

(c) If ALUR points of $B(X^{\ast})$ form a boundary for $X$ (in particular, if
$X^{\ast}$ is ALUR), then $X$ is Fr\'{e}chet smooth.\textquotedblright

\smallskip There is no proof.

\smallskip

\textquotedblleft{Theorem} 6.5 For a Banach space $X$ and $x^{\ast}\in
S(X^{\ast})$, the following are equivalent:

(a) $x^{\ast}$ is a rotund point of $B(X^{\ast})$ as well as a w*-w PC of
$B(X^{\ast})$.

(b) ... (c) ... (d) ... (e) ...

(f) $x^{\ast}$ is a w*-wALUR point of $B(X^{\ast})$. (g) ..."

\smallskip The proof contains 4 $\lim$. The part of the proof referring to (f)
is \textquotedblleft Equivalence of (a) and (f) is immediate from definitions".

\smallskip

\textquotedblleft Replacing the weak topology by the norm topology in the
above Theorem, we immediately obtain

{Corollary} 6.8 Let $X$ be a normed linear space. For $x^{\ast}\in S(X^{\ast
})$, the following are equivalent:

(a) ... (b) ... (c) ... (d) $x^{\ast}$ is a w*-nALUR point of $B(X^{\ast})$.
(e) ..."

\smallskip There is no proof.

\smallskip

\textquotedblleft{Corollary} 6.9 For a Banach space $X$ and $x\in S(X)$, the
following are equivalent:

(a) $x$ is a wALUR point of $B(X)$ as well as a w*-w PC of $B(X^{\ast\ast})$.

(b) $x$ is a wALUR point of $B(X)$ as well as a w*-w seq PC of $B(X^{\ast\ast
})$.

(c) ... (d) ... (e) ... (f) ... (g) ... (h) $x$ is a w*-wALUR point of
$B(X^{\ast\ast})$. (i) ... ."

\smallskip There is no proof.

As a conclusion, in \cite{BaHuLi04}, excepting the definitions,
there is only one instance where $\lim_{m}\lim_{n}$ appears (in the
proof of Prop.\ 2.5), and in all the proofs of the implications
$p\Rightarrow q$ in which ALUR is involved, one uses results from
\cite{BaHuLiTr00} or \cite{BaHuLi04}.

\medskip

-- \cite{AizGar05}: Notice that $\lim$ appears 4 times including $\lim_{m}%
\lim_{n}$ 1 time (on p.\ 132).

\smallskip\textquotedblleft{Remark} 2.9. Recently, Bandyopadhyay and Lin [3]
proved the following equivalences for a point $x$ of the unit sphere of a
Banach space $X$:

(i) ... (ii) $x$ is an almost locally uniformly rotund point of $B_{X}$ if and
only if it is a strongly exposed point of $B_{X}$ for each $f\in S_{X}$ such
that $f(x)=1$."

\smallskip In fact the results mentioned in \cite[Rem.\ 2.9]{AizGar05} are
stated without proofs in \cite{BanLin01}.

\smallskip\textquotedblleft{Corollary} 2.11. Let $X$ be a Banach space. Then
$X$ is locally uniformly rotund if it is almost locally uniformly rotund and
its norm is Fr\'{e}chet differentiable in $S_{X}$."

\smallskip There is no proof.

\medskip

-- \cite{BanGod06}: Notice that $\lim$ appears 3 times including $\lim_{m}%
\lim_{n}$ 1 time (in Definition 2.19).

\smallskip\textquotedblleft Proposition 2.20. Let $X$ be an Asplund space.
Assume that the norm $\left\Vert \cdot\right\Vert _{X}$ on $X$ satisfies one
of the following conditions:

(i) $\left\Vert \cdot\right\Vert _{X}$ is almost LUR. (ii) the dual norm on
$X^{\ast}$ is G\^{a}teaux differentiable.

Then the following are equivalent: (a) $NA(X)$ is spaceable. (b) there exists
a proximinal subspace $Y$ of $X$ such that $X/Y$ is infinite dimensional and reflexive."

\smallskip There is no $\lim$ in the proof. The part of the proof referring to
(i): \textquotedblleft It has been proved in [2] that $x\in S_{X}$ is an
almost LUR point if and only if it is strongly exposed by every functional
that attains its norm at $x$. It follows that (i) holds if and only if ...",
where \textquotedblleft\lbrack2]" is our reference \cite{BaHuLi04}.

\smallskip\textquotedblleft Proposition 2.23. Let $X$ be a Banach space with
the RNP or an almost LUR norm. Then $\operatorname*{span}(NA(X))=X^{\ast}$. ...

Proof. As noted above, if $X$ is almost LUR, then, $NA(X)=...$"

\smallskip There is no $\lim$ in the proof.

\medskip

-- \cite{BaLiLiNa08}: Notice that $\lim$ appears 9 times including $\lim
_{m}\lim_{n}$ 1 time (in Definition 3.1). One mentions \textquotedblleft In
the literature, the acronym ALUR is also used for average locally uniformly
rotund ... In this section, we consider proximinality in $\tau$-ALUR spaces.
... From the results of [1, Proposition 4.4 and Corollary 4.6], it follows
that a Banach space $X$ is $\tau$-ALUR\ $\Leftrightarrow$ every $x^{\ast}\in
NA(X)$ is a $\tau$-strongly exposing functional". Notice that [1] is our
reference \cite{BaHuLi04}.

\smallskip

\textquotedblleft Theorem 3.2. Let $X$ be a $\tau$-ALUR Banach space and $Y$
be a subspace such that $X/Y$ is reflexive. Then the following are equivalent:
(a) ... (b) ... (c) ... (d) ... "

\smallskip There is no $\lim$ in the proof. The part of the proof referring to
$\tau$-ALUR: \textquotedblleft(d)~$\Rightarrow$ (a). Let $Y^{\perp}\subseteq
NA(X).$ Since $X$ is $\tau$-ALUR, every $x^{\ast}\in Y^{\perp}\simeq
(X/Y)^{\ast}$ is in particular an exposing functional."

\smallskip

\textquotedblleft Theorem 3.4. For a Banach space $X$, $X$ is LUR~$\Rightarrow
$ $X$ is ALUR~$\Rightarrow$ $X$ is strictly convex and has the KK property.
And neither converse implication holds."

\smallskip There is no $\lim$ in the proof. The part of the proof referring to
ALUR: \textquotedblleft The fact that LUR~$\Rightarrow$ ALUR, but not
conversely, has been noted in [2]; Assume now that $X$ is ALUR. Then $X$ is
clearly strictly convex". Notice that [2] is our reference \cite{BaHuLiTr00}.

\smallskip\textquotedblleft Remark 3.5. Note that the last part of the proof
shows that if $X$ is not reflexive and an $M$-ideal in $X^{\ast\ast}$, then
$x^{\ast}$ always fails to be wALUR".

\smallskip\textquotedblleft Theorem 3.6. An ALUR Banach space $X$ is
reflexive~$\Leftrightarrow$ the intersection of any two proximinal hyperplanes
is proximinal."

\smallskip There is no $\lim$ in the proof. The part of the proof referring to
ALUR: \textquotedblleft Conversely, let $X$ be ALUR and suppose the
intersection of any two proximinal hyperplanes is proximinal. Let $x^{\ast
},y^{\ast}\in NA(X)$. Then $\ker x^{\ast}$ and $\ker y^{\ast}$ are proximinal
hyperplanes in $X$ ... Since $X$ is ALUR and since the set of strongly
exposing functionals form a $G_{\delta}$ set ...".

\medskip

-- \cite{BaLiRa09}: Notice that $\lim$ appears 6 times including $\lim_{m}%
\lim_{n}$ 1 time (in Definition 2.5(b)).

\smallskip

\textquotedblleft Remark 2.6. ... In view of [1, Corollary 8], in the special
case of $X$ in $X^{\ast\ast}$, we have the following result: If $X$ is wALUR,
in particular, if $X^{\ast}$ is smooth, then no point of $X^{\ast\ast
}\setminus X$ has a farthest point in $B_{X}$." Notice that [1] from the
preceding remark is our reference \cite{BaHuLiTr00}.

\smallskip There is no proof.

\medskip

-- \cite{Pau10}: Notice that $\lim$ appears 16 times including $\lim_{m}%
\lim_{n}$ 1 time (in Definition 6.6(c)).

\smallskip\textquotedblleft Theorem 7.2.12. $X$ is wALUR if and only if $X$ is
rotund and ARB.

Proof. We recall [7, Corollary 8] that $x\in S_{X}$ is a wALUR point
of $B_{X}$ if and only if $x$ is a rotund point of
$B_{X^{\ast\ast}}$. Let $X$ be wALUR. Then $X$ is clearly rotund.
..."

\smallskip There is no $\lim$ in the proof. Notice that [7] is our reference
\cite{BaHuLiTr00}.

\medskip

-- \cite{BanPau10}: Notice that $\lim$ appears 6 times including $\lim_{m}%
\lim_{n}$ 1 time (in Definition 6.6(c)).

\smallskip\textquotedblleft Theorem 6.7. $X$ is wALUR if and only if $X$ is
rotund and ARB.

Proof. We recall [2, Corollary 8] that $x\in S_{X}$ is a wALUR point of
$B_{X}$ if and only if $x$ is a rotund point of $B_{X^{\ast\ast}}$. Let $X$ be
wALUR. Then $X$ is clearly rotund. ..."

\smallskip There is no $\lim$ in the proof. Notice that [2] is our reference
\cite{BaHuLiTr00}.

\medskip

-- \cite{ZhaLiu11}: Notice that $\lim$ appears 37 times including $\lim
_{m}\lim_{n}$ 1 time (on p. 3897).

\smallskip\textquotedblleft Theorem 3.6. A Banach space $X$ is very convex if
and only if every $x^{\ast}\in D(S(X))$ is a smooth point of $B(X^{\ast})$."

\smallskip There are 2 $\lim$ in the proof. The part of the proof referring to
ALUR: \textquotedblleft Bandgapadhyay et al.\ (Proposition 4.4 in [11]) proved
that: $x$ is a wALUR point of $B(X)$ if and only if every $x^{\ast}\in D(x)$
is a smooth point of $B(X^{\ast})$. Hence we may know that : $X$ is wALUR if
and only if for any $x^{\ast}\in D(S(X))$ is a smooth point of $B(X^{\ast})$."
Notice that [11] is our reference \cite{BaHuLi04}.

\smallskip\textquotedblleft Theorem 3.7. A Banach space $X$ is very convex if
and only if it is wALUR."

\smallskip There is no proof.

\smallskip\textquotedblleft Theorem 3.16. A Banach space $X$ is strongly
convex if and only if it is ALUR."

\smallskip Before the statement of \cite[Theorem 3.16]{ZhaLiu11} one finds the
following text in which one refers to ALUR: \textquotedblleft Bandgapadhyay et
al.\ (Corollary 4.6 in [11]) proved that: $x$ is a wALUR point of $B(X)$ if
and only if every $x^{\ast}\in D(x)$ is a Fr\'{e}chet smooth point of
$B(X^{\ast})$. Hence we may know that: $X$ is ALUR if and only if for any
$x^{\ast}\in D(S(X))$ is a Fr\'{e}chet smooth point of $B(X^{\ast})$. By this
conclusion and Corollary 3.15, we can obtain the following theorem."

\medskip

-- \cite{ZhaLiu11b}: Notice that $\lim$ appears 9 times including $\lim
_{m}\lim_{n}$ 1 time (on p. 2).

\smallskip There are no results referring to ALUR. However, one finds the
following text on page 2: \textquotedblleft Recently, we have proved that ALUR
and strong convexity, WALUR and very convex are equivalent, respectively
[8].", where [8] is our reference \cite{ZhaLiu11}.

\medskip

-- \cite{ZhaLiu12}: Notice that $\lim$ appears 10 times including $\lim
_{m}\lim_{n}$ 1 time (on p. 2).

\smallskip There are no results referring to ALUR. However, one finds the
following text on page 2: \textquotedblleft Recently, we proved that almost
locally uniformly rotund space is equivalent to strongly convex space and that
weakly almost locally uniformly rotund space is equivalent to very convex
space [7].", where [7] is our reference \cite{ZhaLiu11}.

\medskip

-- \cite{DamBaj17}: Notice that $\lim$ appears 23 times including $\lim
_{m}\lim_{n}$ 1 time (Definition 6.1).

\smallskip\textquotedblleft Remarks 6.2. We have implication LUR $\Rightarrow$
ALUR $\Rightarrow$ R (rotund)."

\smallskip There is no proof.

\smallskip\textquotedblleft Remarks 6.3. If $x\in S(X)$ is a strongly exposed
point and a smooth point of $B(X)$, then it is an almost locally uniformly
rotund point of $B(X)$."

\smallskip There is no proof.

\smallskip\textquotedblleft Theorem 6.5. Aizpurlt and Garcia-Pacheco (2005).
Let $X$ be a Banach space. Then $X$ is locally uniformly rotund if
it is almost locally uniformly rotund and its norm is Frechet
differentiable in $S(X)$."

\smallskip There is no $\lim$ in the proof. The part of the proof referring to
ALUR: `Assume that a Banach space $X$ is ALUR and its norm is Frechet
differentiable. ... Then $x$ is strongly exposed point of $B(X)$ because
\textquotedblleft$x$ is an almost locally uniformly rotund point of $B(X)$ if
and only if it is a strongly exposed point of $B(X)$ for each $f\in S(X^{\ast
})$ such that $f(x)=1$\textquotedblright. ... Thus, by Theorem 6.4 norm of $X$
is LUR.'

\medskip

-- \cite{GupNar17}: Notice that $\lim$ appears 13 times including $\lim
_{m}\lim_{n}$ 1 time (on p.\ 114).

\smallskip\textquotedblleft Theorem 3.4. Let $X$ be Banach space. Then the
following are equivalent:

(i) ... (ii) Every $x\in NA(X)$ is strongly exposing functional. (iii) ...

(iv) $X$ is almost locally uniform rotund."

\smallskip There is no $\lim$ in the proof. The part of the proof referring to
ALUR: \textquotedblleft(ii)~$\Leftrightarrow$ (iv) follows from Corollary 4.6
of [2]". Notice that [2] is our reference \cite{BaHuLi04}.

\medskip

-- \cite{JuMaRu23}: Notice that $\lim$ appears 7 times including $\lim_{m}%
\lim_{n}$ 1 time (on p. 20).

\smallskip\textquotedblleft Corollary 3.12. ALUR Banach spaces satisfy
property [P]."

\smallskip There is no proof. However, before Corollary 3.12, one says
\textquotedblleft It is clear that LUR spaces are ALUR, but the reverse
implication is not true (see [16, Corollary 12]). It is observed in [15,
Corollary 4.6] that if $X$ is ALUR, then each point $x$ in $S_{X}$ is strongly
exposed by every $x^{\ast}\in S_{X^{\ast}}$ which attains its norm at $x$.
Thus, in particular, if $X$ is ALUR, then $S_{X}={}$str-exp($B_{X}$)". Notice
that [16] and [15] are our references \cite{BaHuLiTr00} and \cite{BaHuLi04}, respectively.

\medskip

-- \cite{DeBSom24}: Notice that $\lim$ appears 16 times including $\lim
_{m}\lim_{n}$ 1 time (on p.\ 1692); also notice that ALUR is the abbreviation
for \textquotedblleft average LUR" in \cite{DeBSom24}.

\smallskip On page 1689 one says \textquotedblleft It is well-known and
easy-to-prove that local uniform rotundity of the unit ball implies that each
point of unit sphere is strongly exposed by any of its supporting functionals
(spaces satisfying this condition are called almost locally uniformly rotund
(almost LUR) in the literature [1, 2])". Notice that [1, 2] are our references
\cite{AizGar05, BanLin01}, respectively.

\smallskip\textquotedblleft Theorem 2.4. Let $X$ be a Banach space, then the
following holds true

(i) if $X$ is almost LUR, then $X$ is ALUR.

(ii) if $X$ is reflexive and ALUR, then $X$ is almost LUR."

\smallskip There is no $\lim$ in the proof. The parts of the proof referring
to almost LUR: \textquotedblleft(i) It follows by the fact that if $x\in
S_{X}$ is strongly exposed by $x^{\ast}\in X^{\ast}$, then $x$ is a denting
point of $B_{X}$"; \textquotedblleft(ii) By the \v{S}mulyan Lemma, $x$ is
strongly exposed by $x^{\ast}$, hence $X$ is almost LUR."

\smallskip\textquotedblleft Example 2.5. Let $X=\ell_{1}$ endowed with the
norm $\left\Vert (x_{n})_{n}\right\Vert =...$ Let us show that $\left\Vert
\cdot\right\Vert $ is not almost LUR. ... It remains to show that $x$ is not
strongly exposed by $x^{\ast}.$ ... Observing that $\left\Vert x-\left(
\frac{n}{1+n}\right)  e_{n}\right\Vert >1/2$, we get that $x$ is not strongly
exposed by $x^{\ast}$, hence $\left\Vert \cdot\right\Vert $ is not almost LUR".

\medskip

-- \cite{DasPau24}: Notice that $\lim$ appears 19 times including one
$\lim_{m,n}$ (in Definition 1.10(b)) and one $\lim_{n,m}$ (in the proof of
Theorem 2.16); observe that there is no a definition of $\lim_{m,n}%
\gamma_{m,n}$ or $\lim_{n,m}\gamma_{m,n}$ (with $\gamma_{m,n}=h_{m}\left(
\frac{x_{n}+x}{2}\right)  $ in both cases). Having in view that one refers to
\cite{BaLiLiNa08} in \cite[Def.\ 1.10(b)]{DasPau24} (and one refers to
\cite{BaHuLiTr00} in \cite[Def.\ 3.1]{BaLiLiNa08}), one could consider that
$\lim_{m,n}h_{m}\left(  \frac{x_{n}+x}{2}\right)  $ coincides with $\lim
_{m}\lim_{n}h_{m}\left(  \frac{x_{n}+x}{2}\right)  $ in \cite{DasPau24}, as in
\cite[Def.\ 2(c)]{BaHuLiTr00} and \cite[Def.\ 3.1]{BaLiLiNa08}); even more, in
these references (and \cite{BaHuLi04}) $\{x_{m}^{\ast}\}\subseteq B_{X^{\ast}%
}$, and so, if one wishes to apply results from such references (see Remark
2.5) one must pay attention because $(h_{m})\subseteq S_{X^{\ast}}$!

\smallskip\textquotedblleft Remark 2.5 In a $\tau$-ALUR space a norm attaining
functional is a $\tau$-strongly exposing functional (see [19, Proposition 4.4,
Corollary 4.6]). By virtue of [11, Theorem 2.11] the condition stated in
Theorem 2.2 (b) follows if $X$ is ALUR and $f\in NA(X)$. Hence a Banach space
that is reflexive and ALUR is an example where every hyperplane (and hence
every closed convex set) is approximatively compact. Our next result is a
stronger version of this observation."

Notice that [19] and [11] from the quoted text are our references
\cite{BaHuLi04} and \cite{BaLiLiNa08}, respectively.

\smallskip\textquotedblleft Theorem 2.16 Let $X$ be a Banach space that is
wALUR(wLUR), then The following are equivalent.

(a) Every proximinal hyperplane of $X$ is strongly proximinal.

(b) Every proximinal hyperplane of $X$ is approximatively compact.

(c) $X$ is ALUR(LUR).

Proof We derive the case for wALUR spaces; similar arguments work for wLUR
spaces. From Remark 2.5 we have (c)~$\Rightarrow$ (b) and (b)~$\Rightarrow$
(a). We prove (a)~$\Rightarrow$ (c). Let us assume that $x\in S_{X}$,
$(x_{n})\subseteq B_{X}$ and $(h_{m})\subseteq S_{X^{\ast}}$ be such that
$\lim_{n,m}h_{m}\left(  \frac{x_{n}+x}{2}\right)  =1$. Let $f\in S_{X^{\ast}}$
be such that $f(x)=1$. As $X$ is wALUR, we obtain $|f(x_{n})-f(x)|\rightarrow
0$. From our assumption (a), it follows that $f$ is an SSD point of $X^{\ast}%
$, and this leads to $d(x_{n},J_{X}(f))\rightarrow0$. As $X$ is strictly
convex we have $\left\Vert x_{n}-x\right\Vert \rightarrow0$. This completes
the proof. $\square$"

\smallskip Having in view Definition 1.10(b) and the aim of proving (c) when
(a) holds, one may consider that \textquotedblleft$\lim_{n,m}$" instead of
\textquotedblleft$\lim_{m,n}$" is a misprint in the above proof (what we do).
So, in order to prove (a)~$\Rightarrow$ (c) one considers \textquotedblleft%
$x\in S_{X}$, $(x_{n})\subseteq B_{X}$ and $(h_{m})\subseteq
S_{X^{\ast}}$ such that $\lim_{m,n}h_{m}\left(
\frac{x_{n}+x}{2}\right)  =1$". Because $X$ is wALUR one has
$x_{n}\rightarrow^{w}x$ by Definition 1.10(b), and one must prove
that $x_{n}\rightarrow x$. For this one takes \textquotedblleft$f\in
S_{X^{\ast}}$ ... such that $f(x)=1$" (and so $f\in NA(X)$). Is
$\ker f$ proximinal?\footnote{Thanks to Prof.\ P. Shunmugaraj for
letting us know the following result: Let $f\in S_{X^{\ast}};$ then
$f\in NA(X)$ if and only if $\ker f$ is proximinal.}

\medskip

-- \cite{ZhZhLi17}: Notice that $\lim$ appears 29 times including $\lim
_{m}\lim_{n}$ 9 times (two of them in Definition 1.1). Observe that one refers
to our reference \cite{BaHuLiTr00} for the definition of WALUR point in
\cite[Def.\ 1.1(2)]{ZhZhLi17}, but WALUR is not mentioned later on. In
\cite[Def.\ 1.1(2$^{\prime}$)]{ZhZhLi17} one says that $x\in S(X)$ is
\textquotedblleft a CWALUR point of $B(X)$ if, for any $\{x_{n}\}_{n=1}%
^{\infty}\subseteq B(X)$ and $\{x_{m}^{\ast}\}_{m=1}^{\infty}\subseteq
B(X^{\ast})$, the condition $\lim_{m}\lim_{n}x_{m}^{\ast}\left(  \frac
{x+x_{n}}{2}\right)  =1$ implies that $\{x_{n}\}_{n=1}^{\infty}$ is relatively
weakly compact."

\smallskip Besides the meaning of $\lim_{m}\lim_{n}\gamma_{m,n}$, in the
sequel we are also interested to understand what is meant by \textquotedblleft%
$\{x_{n}\}_{n=1}^{\infty}$ is relatively weakly compact".

Notice that in Definition 1.1 (2$^{^{\prime}}$) and \textquotedblleft(3) (see
[11], [12])" one speaks about the sequence $\{x_{n}\}_{n=1}^{\infty}$ as being
\textquotedblleft relatively weakly compact" or \textquotedblleft relatively
compact", while in \textquotedblleft(4) (see [9], [10])" one speaks about the
sequence $\{x_{n}^{\ast}\}_{n=1}^{\infty}$ as being \textquotedblleft
relatively weakly compact" or \textquotedblleft relatively compact". So,
\textit{which is the meaning of (relatively) ($\tau$-)compact sequence?}

The mentioned references are:

``9. B.\ B.\ Panda and O.\ P.\ Kapoor, A generalization of local
uniform convexity of the norm, J.\ Math.\ Anal.\ Appl.\ 52 (1975),
300-305. Zbl 0314.46014.\ MR0380365.

10. J.\ Wang and C.\ Nan, On the dual spaces of the S-spaces and WkR spaces,
Chinese J. Contemp.\ Math.\ 13 (1992), no.\ 1, 23-27. Zbl 0779.46023. MR1239303.

11. J.\ Wang and Z.\ Zhang, Characterizations of the property (C-$\kappa$),
Acta Math.\ Sci.\ Ser.\ A Chin.\ Ed.\ 17 (1997), no.\ 3, 280-284. Zbl
0917.46013. MR1484502.

12. Z.\ Zhang and Z.\ Shi, Convexities and approximative compactness
and continunity of metric projection in Banach spaces, J.\ Approx.\
Theory. 161 (2009), no.\ 2, 802-812. Zbl 1190.46018. DOI
10.1016/j.jat.2009.01.003."

\medskip

We succeeded to find Ref.\ 9 and Ref.\ 12.

\smallskip

In Ref.\ 9 the string \textquotedblleft relative" is absent. In Theorem 2.2 of
Ref.\ 9 one finds the text \textquotedblleft Suppose that $\{g_{n}\}$ is a
sequence in $X$ with ... Then the sequence $\{g_{n}\}$ has a convergent
subsequence.", while in its proof one finds \textquotedblleft As $X$ has the
property (M), the sequence $\{z_{n}\}$ has a convergent subsequence. Clearly,
this establishes the compactness of the sequence $\{g_{n}\}$ and the proof is complete."

In Ref.\ 12 one finds \textquotedblleft A Banach space $X$ is said to be ...
nearly strongly convex ... [11,18] if for any $x\in S(X)$ and $\{x_{n}%
\}_{n\in\mathbb{N}}\subset B(X)$ ... then $\{x_{n}\}_{n\in\mathbb{N}}$ is
relatively compact"; \textquotedblleft Lemma 2.1. Suppose that $X$ is a nearly
strongly convex Banach space ... then $\{y_{n}\}_{n\in\mathbb{N}}$ is
relatively compact and each of its convergent subsequences converges to an
element of $P_{A}(x)$. ... We will prove ... that $\{y_{n}\}_{n\in\mathbb{N}}$
is relatively compact. ... Since $X$ is nearly strongly convex, we know that
$\{-\frac{y_{n}}{\left\Vert y_{n}\right\Vert }\}$ is a \emph{relatively
compact set} ... we obtain that $\{y_{n}\}_{n\in\mathbb{N}}$ i\emph{s
relatively compact}."

Ref.\ 12 is also cited in \cite{ZhZhLi16} written by the authors of
\cite{ZhZhLi17} (the paper in the present discussion). In \cite[page
603]{ZhZhLi16} one finds \textquotedblleft... $\{\frac{x_{n}}{\left\Vert
x_{n}\right\Vert }\}_{n=1}^{\infty}$ and therefore $\{x_{n}\}_{n=1}^{\infty}$
\emph{are relatively weakly compact set}. It follows that $\{x_{n}%
\}_{n=1}^{\infty}$ has a weakly convergent subsequence".\footnote{On
05.02.2025 we sent a message to the authors of \cite{ZhZhLi17} containing the
following two questions related to \cite[Def.\ 1.1]{ZhZhLi17}: 1) How is
defined the iterated limit in (2) and (2$^{\prime}$)? Please use the
quantifiers $\forall$ and $\exists$. 2) What do you mean by the fact that the
sequence $\{x_{n}\}$ is relatively compact with respect to a topology (for
example, weak, weak* or strong)? No answer until now.}

\smallskip Even if it is not clear from Ref.\ 9, Ref.\ 12 and \cite{ZhZhLi16}
what is meant by a relatively compact sequence, in the sequel we consider that
a sequence $\{y_{n}\}$ from a topological space $(Y,\tau)$ is relatively
($\tau$-)~compact if $\{y_{n}\}$ has (at least) a ($\tau$-)~converging
subsequence $\{y_{n_{k}}\}$. (However, I should say that $\{y_{n}\}$ is
relatively ($\tau$-)~compact if each subsequence of $\{y_{n}\}$ has, at its
turn, a convergent subsequence; otherwise, the sequence $\{y_{n}%
\}\subset(\mathbb{R},\tau_{0})$ with $y_{n}:=n$ for $n$ even and $y_{n}:=1/n$
for $n$ odd would be relatively compact!)

\medskip

\textquotedblleft Theorem 2.1. Let $X$ be a Banach space, with $x\in S(X)$.
Then the following are equivalent:

(1) $x$ is a nearly rotund point of $B(X^{\ast\ast})$;

(2) for any $\{x_{n}^{\ast\ast}\}_{n=1}^{\infty}\subseteq B(X^{\ast\ast})$ and
$\{x_{m}^{\ast}\}_{m=1}^{\infty}\subseteq B(X^{\ast})$, the condition
$\lim_{m}\lim_{n}\left(  \frac{x+x_{n}^{\ast\ast}}{2}\right)  (x_{m}^{\ast
})=1$ implies that all w$^{\ast}$-cluster points of $\{x_{n}^{\ast\ast
}\}_{n=1}^{\infty}$ belong to $X$;

(3) for every unbounded nested sequence $\{B_{n}^{\ast}\}_{n=1}^{\infty}$ of
balls such that $x$ is bounded below on $\cup B_{n}^{\ast}$, if for any
$\{y_{n}^{\ast\ast}\}_{n=1}^{\infty}\subset S(X^{\ast\ast})$ the sequence
$\{\inf y_{n}^{\ast\ast}(B_{n}^{\ast})\}_{n=1}^{\infty}$ is also bounded
below, then all w$^{\ast}$-cluster points of $\{y_{n}^{\ast\ast}%
\}_{n=1}^{\infty}$ belong to $X$;

(4) for every unbounded nested sequence $\{B_{n}^{\ast}\}_{n=1}^{\infty}$ of
balls in $X^{\ast}$ such that $x$ is bounded below on $\cup B_{n}^{\ast}$, if
for any $\{y_{n}\}_{n=1}^{\infty}\subset S(X)$ the sequence $\{\inf
y_{n}(B_{n}^{\ast})\}_{n=1}^{\infty}$ is also bounded below, then
$\{y_{n}\}_{n=1}^{\infty}$ is relatively weakly compact;

(5) for every unbounded nested sequence $\{B_{n}^{\ast}\}_{n=1}^{\infty}$ of
balls in $X^{\ast}$ such that $x$ is bounded below on $\cup B_{n}^{\ast}$, if
$x^{\ast\ast}\in S(X^{\ast\ast})$ is also bounded below on $\cup B_{n}^{\ast}$
then $x^{\ast\ast}\in X$.

\smallskip Proof. The proof is similar to Theorem 6 in [2]. We will give their
proof for the sake of completeness.

(1)~$\Rightarrow$ (2) Let $\{x_{n}^{\ast\ast}\}_{n=1}^{\infty}\subseteq
B(X^{\ast\ast})$, $\{x_{m}^{\ast}\}_{n=1}^{\infty}\subseteq B(X^{\ast})$ such
that $\lim_{m}\lim_{n}\left(  \frac{x+x_{n}^{\ast\ast}}{2}\right)
(x_{m}^{\ast})=1$. Let $y^{\ast\ast}$ be a w*-cluster point of $\{x_{n}%
^{\ast\ast}\}_{n=1}^{\infty}$; then $\lim_{m}\left(  \frac{x+y^{\ast\ast}}%
{2}\right)  (x_{m}^{\ast})=1$, and hence, $\left\Vert \frac{x+y^{\ast\ast}}%
{2}\right\Vert =1$. By (1), $y^{\ast\ast}\in X$.

(2)~$\Rightarrow$ (3) ... The fact that $\{B_{n}^{\ast}\}$ is nested implies
that $y_{n}^{\ast\ast}(y_{m}^{\ast})\geq1+\frac{c}{r_{m}}$ for all $n\geq m$.
It follows that $\left(  \frac{x+y_{n}^{\ast\ast}}{2}\right)  (y_{m}^{\ast
})\geq1+\frac{c}{r_{m}}$ for all $n\geq m$. Since $r_{m}\rightarrow\infty$, we
conclude with $\lim_{m}\lim_{n}\left(  \frac{x+y_{n}^{\ast\ast}}{2}\right)
(y_{m}^{\ast})=1$. By (2), we obtain that all w*-cluster points of
$\{y_{n}^{\ast\ast}\}_{n=1}^{\infty}$ belong to $X$.

(3)~$\Rightarrow$ (4) This is clear.

(4)~$\Rightarrow$ (5) ... By (4), we know that $\{y_{n}\}_{n=1}^{\infty}$ is
relatively weakly compact. Therefore, there exists a subsequence $\{y_{n_{k}%
}\}_{k=1}^{\infty}\subseteq\{y_{n}\}_{n=1}^{\infty}$ and $y\in X$ such that
$y_{n_{k}}\rightarrow^{w}y$ and so we may assure that $x_{n_{k}}$ is also
weakly convergent."

\medskip Before discussing the above text, recall that, by Definition
1.1.(1$^{^{\prime}}$), \textquotedblleft$x\in S(X)$ is a nearly rotund point
of $B(X^{\ast\ast})$ if $x^{\ast\ast}\in X$ for any $x^{\ast\ast}\in
X^{\ast\ast}$ satisfying $\left\Vert x^{\ast\ast}\right\Vert =\left\Vert
\frac{x+x^{\ast\ast}}{2}\right\Vert =1$".

\begin{remark}
\label{R1} \emph{(i)} Observe that the proofs of the implications
(1)~$\Rightarrow$ (2) and (2)~$\Rightarrow$ (3) mimic those of the
implications (a)~$\Rightarrow$ (b) and (b)~$\Rightarrow$ (c) of \cite[Th.\ 6]%
{BaHuLiTr00}, respectively, and so our remarks are the same. Moreover, the
proofs of (4)~$\Rightarrow$ (5) and (5)~$\Rightarrow$ (1) mimic those of the
implications (c$^{\prime}$)~$\Rightarrow$ (d) of \cite[Cor.\ 8]{BaHuLiTr00}
and (e)~$\Rightarrow$ (a) of \cite[Th.\ 6]{BaHuLiTr00}, respectively.

\emph{(ii)} One says: \textquotedblleft(3)~$\Rightarrow$ (4) This is
clear." Although let us try a proof. Set $\widehat{u}:=F_{X}(u)$ for
$u\in X$, where $F_{X}:X\rightarrow X^{\ast\ast}$ is the canonical
injection. Consider $y_{n}^{\ast\ast}:=\widehat{y_{n}}$ $(\in
B_{X^{\ast\ast}})$ and take
$\{y_{n_{i}}^{\ast\ast}\}_{i\in I}$ a subnet of $\{y_{n}^{\ast\ast}%
\}_{n=1}^{\infty}$. Because $B_{X^{\ast\ast}}$ is w*-compact, $\{y_{n_{i}%
}^{\ast\ast}\}$ has a subnet $\{y_{n_{i_{j}}}^{\ast\ast}\}_{j\in J}$
w*-converging to $y^{\ast\ast}$ $(\in B_{X^{\ast\ast}})$; hence $y^{\ast\ast}$
is a cluster point of $\{y_{n}^{\ast\ast}\}_{n=1}^{\infty}$, and so
$y^{\ast\ast}=\widehat{y}$ for some $y\in B_{X}$ by (3). Clearly, $y$ is a
w-cluster point of $\{y_{n}\}$. Is $y$ the w-limit of a subsequence of
$\{y_{n}\}_{n=1}^{\infty}$?

\emph{(iii)} In the proof of (4)~$\Rightarrow$ (5) one obtained:
\textquotedblleft Therefore, there exists a subsequence $\{y_{n_{k}}%
\}_{k=1}^{\infty}\subseteq\{y_{n}\}_{n=1}^{\infty}$ and $y\in X$ such that
$y_{n_{k}}\rightarrow^{w}y$" because \textquotedblleft$\{y_{n}\}_{n=1}%
^{\infty}$ is relatively weakly compact".
\end{remark}

\textquotedblleft Theorem 2.2. Let $X$ be a Banach space, with $x\in S(X)$.
The following are equivalent:

(1) $x$ is a nearly rotund point of $B(X^{\ast\ast})$;

(2$^{\prime}$) $x$ is a CWALUR point of $B(X)$;

(3$^{\prime}$) for every $x^{\ast\ast}\in S(X^{\ast\ast})$, if there exists
sequence $\{x_{n}^{\ast}\}_{n=1}^{\infty}\subseteq S(X^{\ast})$ such that
$\lim_{n}x^{\ast\ast}(x_{n}^{\ast})=\lim_{n}x_{n}^{\ast}(x)=1$, then
$x^{\ast\ast}\in X$;

(4$^{\prime}$) for any $\{x_{n}^{\ast\ast}\}_{n=1}^{\infty}\subseteq
B(X^{\ast\ast})$, if $\left(  \frac{x+x_{n}^{\ast\ast}}{2}\right)
_{n=1}^{\infty}$ is asymptotically normed by $B(X^{\ast})$, then all $w^{\ast
}$-cluster points of $\{x_{n}^{\ast\ast}\}_{n=1}^{\infty}$ belong to $X$.

\smallskip

Proof. (1)~$\Leftrightarrow$ (2$^{\prime}$) Let $\{x_{n}\}_{n=1}^{\infty
}\subseteq B(X)$ and $\{x_{m}^{\ast}\}_{m=1}^{\infty}\subseteq B(X^{\ast})$
such that $\lim_{m}\lim_{n}x_{m}^{\ast}\left(  \frac{x+x_{n}}{2}\right)  =1$.
By the relatively weak$^{\ast}$-compactness of $\{x_{n}\}_{n=1}^{\infty
}\subseteq B(X^{\ast\ast})$, we obtain a w$^{\ast}$-cluster point $x^{\ast
\ast}$ of $\{x_{n}\}_{n=1}^{\infty}$. Then, it is easy to prove that $\lim
_{m}\lim_{n}x_{m}^{\ast}\left(  \frac{x+x^{\ast\ast}}{2}\right)  =1$, which
further implies that $x^{\ast\ast}\in X$ by (1). Consequently, $\{x_{n}%
\}_{n=1}^{\infty}$ is relatively weakly compact. This shows $x$ is a CWALUR
point of $B(X)$."

\smallskip Let us have a closer look at the proof of (1)~$\Leftrightarrow$
(2$^{\prime}$) above (in fact, of (1)~$\Rightarrow$ (2$^{\prime}$)).

Condition (1) of Theorem 2.2 coincides with condition (1) of Theorem 2.1.
Taking $x_{n}^{\ast\ast}:=\widehat{x_{n}}$ one obtains, as in the proof of the
implication (1)~$\Rightarrow$ (2) of Theorem 2.1, \textquotedblleft a
w$^{\ast}$-cluster point $y^{\ast\ast}$ of $\{x_{n}\}_{n=1}^{\infty}$; then
$\lim_{m}\left(  \frac{x+y^{\ast\ast}}{2}\right)  (x_{m}^{\ast})=1$, and
hence, $\left\Vert \frac{x+y^{\ast\ast}}{2}\right\Vert =1$. By (1),
$y^{\ast\ast}\in X$."

Let us set $x^{\ast\ast}:=y^{\ast\ast}$. So, it is not necessary to mention
\textquotedblleft it is easy to prove that $\lim_{m}\lim_{n}x_{m}^{\ast
}\left(  \frac{x+x^{\ast\ast}}{2}\right)  =1$" but directly $\lim_{m}\left(
\frac{x+x^{\ast\ast}}{2}\right)  (x_{m}^{\ast})=1$ instead of $\lim_{m}%
\lim_{n}x_{m}^{\ast}\left(  \frac{x+x^{\ast\ast}}{2}\right)  =1,$ which does
not look very nice because a) $x_{m}^{\ast}\left(  \frac{x+x^{\ast\ast}}%
{2}\right)  $ makes sense only if one identifies $X^{\ast}$ with $F_{X^{\ast}%
}(X^{\ast})$ (otherwise it must be written $\left(  \frac{x+x^{\ast\ast}}%
{2}\right)  (x_{m}^{\ast})$) and b) $x_{m}^{\ast}\left(  \frac{x+x^{\ast\ast}%
}{2}\right)  $ does not depend on $n$ and so $\lim_{n}x_{m}^{\ast}\left(
\frac{x+x^{\ast\ast}}{2}\right)  =x_{m}^{\ast}\left(  \frac{x+x^{\ast\ast}}%
{2}\right)  $!

\begin{remark}
\label{R2} \emph{(i)} Having $X$ a Banach space, $x\in S(X)$ a nearly rotund
point of $B(X^{\ast\ast})$, $\{x_{n}\}_{n=1}^{\infty}\subseteq B(X)$ and
$\{x_{m}^{\ast}\}_{m=1}^{\infty}\subseteq B(X^{\ast})$ such that $\lim_{m}%
\lim_{n}x_{m}^{\ast}\left(  \frac{x+x_{n}}{2}\right)  =1$, is $\{x_{n}%
\}_{n=1}^{\infty}$ relatively weak$^{\ast}$-compact, that is, there exists a
subsequence $\{x_{n_{k}}\}_{k=1}^{\infty}$ of $\{x_{n}\}_{n=1}^{\infty}$
w*-converging to some $x^{\ast\ast}\in X^{\ast\ast}$?

\emph{(ii)} The proof of (2$^{\prime}$)~$\Rightarrow$ (1) mimics that of
\cite[Th.\ 6, (b)~$\Rightarrow$ (c)]{BaHuLiTr00}, while the proof of
(1)~$\Leftrightarrow$ (4$^{\prime}$) mimics that of \cite[Prop.\ 2.5,
(a)~$\Leftrightarrow$ (b)]{BaHuLi04}
\end{remark}

\textbf{Questions}. (i) \emph{Which is the definition of }$\lim_{m}\lim
_{n}\gamma_{m,n}=\gamma$\emph{ for }$\gamma,$\emph{ }$\gamma_{m,n}%
\in\mathbb{R}$ \emph{in the works dealing with almost LUR (ALUR) or CWALUR?}

\smallskip(ii) \emph{Which is the definition of a relatively compact sequence
in a topological space?}

\medskip

The previous discussion on the results concerning almost LUR points or spaces,
respectively CWALUR points, shows that practically there are no proofs that
directly use the definitions of such points (spaces) excepting those from
\cite{BaHuLiTr00} and \cite{BaHuLi04}. Even in these papers $\lim_{m}\lim_{n}$
does appear only in the proofs of \cite[Th.\ 6]{BaHuLiTr00} and
\cite[Prop.\ 2.5]{BaHuLi04}.

\medskip We propose to replace $\lim_{m}\lim_{n}$ by $\lim_{m}\liminf_{n}$ in
all definitions of almost LUR points (and the similar ones). For example,
\cite[Defs.\ 2(c) and 5]{BaHuLiTr00} become:

\medskip\textbf{Definition B}. Let $X$ be a Banach space. We say that $x\in
S_{X}$ is:

\smallskip(i) an almost LUR (ALUR) (resp.\ weakly almost LUR (WALUR)) point of
$B_{X}$ if for any $\{x_{n}\}\subseteq B_{X}$ and $\{x_{m}^{\ast}\}\subseteq
B_{X^{\ast}}$, the condition
\[
\lim_{m\rightarrow\infty}\liminf_{n\rightarrow\infty}x_{m}^{\ast}\left(
\frac{x_{n}+x}{2}\right)  =1
\]
implies $\lim_{n}\left\Vert x_{n}-x\right\Vert =0$ (resp.\ w-$\lim_{n}%
(x_{n}-x)=0$);

\smallskip(ii) a w*-ALUR point of $B_{X^{\ast}}$ if for any $\{x_{n}^{\ast
}\}\subseteq B_{X^{\ast}}$ and $\{x_{m}\}\subseteq B_{X}$ the condition{%
\begin{equation}
\lim_{m\rightarrow\infty}\liminf_{n\rightarrow\infty}\left(  \frac{x_{n}%
^{\ast}+x^{\ast}}{2}\right)  (x_{m})=1 \label{rz1b}%
\end{equation}
} implies w*-$\lim x_{n}^{\ast}=x^{\ast}$.

\medskip

The following definitions are \cite[Def.\ 1]{BaHuLiTr00} and \cite[Def.\ 1.8
(a), (b)]{BaHuLi04}, respectively.

\medskip\textbf{Definition C}. An \textquotedblleft unbounded nested sequence
of balls\textquotedblright\ in a Banach space $X$ is an increasing sequence
$\{B_{n}=B(x_{n},r_{n})\}$ of open balls in $X$ with $r_{n}\uparrow\infty$.

\medskip\textbf{Definition D.} (a) A subset $\Phi$ of $B_{X^{\ast}}$ is called
a norming set for $X$ if $\left\Vert x\right\Vert =\sup_{x^{\ast}\in\Phi
}x^{\ast}(x)$, for all $x\in X$.

\smallskip(b) A sequence $\{x_{n}\}$ in $B_{X}$ is said to be asymptotically
normed by $\Phi\subseteq B_{X^{\ast}}$ if for any $\varepsilon>0$, there
exists a $x^{\ast}\in\Phi$ and $N\in\mathbb{N}$ such that $x^{\ast}%
(x_{n})>1-\varepsilon$ for all $n\geq N$.

\medskip

In the next result one uses Definitions B, C and D. In this result condition
(b$^{\prime}$) is condition (b) from \cite[Prop.\ 2.5]{BaHuLi04}; removing
(b$^{\prime}$) remains the statement of \cite[Th.\ 6]{BaHuLiTr00}.

\medskip\textbf{Theorem E} (\cite[Th.\ 6]{BaHuLiTr00} \& \cite[Prop.\ 2.5]%
{BaHuLi04}). Let $X$ be a Banach space. For $x^{\ast}\in S_{X^{\ast}}$, the
following are equivalent:

\smallskip(a) $x^{\ast}$ is a rotund point of $B_{X^{\ast}}$;

\smallskip(b) $x^{\ast}$ is a w*-ALUR point of $B_{X^{\ast}}$;

\smallskip(b$^{\prime}$) For any $\{x_{n}^{\ast}\}\subseteq B_{X^{\ast}}$, if
$\{(x_{n}^{\ast}+x^{\ast})/2\}$ is asymptotically normed by $B_{X}$, then
w*-$\lim x_{n}^{\ast}=x^{\ast}$;

\smallskip(c) for every unbounded nested sequence $\{B_{n}\}$ of balls such
that $x^{\ast}$ is bounded below on $\cup B_{n}$, if for any $\{y_{n}^{\ast
}\}\subseteq S_{X^{\ast}}$, the sequence $\{\inf y_{n}^{\ast}(B_{n})\}$ is
bounded below, then w*-$\lim y_{n}^{\ast}=x^{\ast}$;

\smallskip(d) for every unbounded nested sequence $\{B_{n}\}$ of balls such
that $x^{\ast}$ is bounded below on $\cup B_{n}$, if $y^{\ast}\in S_{X^{\ast}%
}$ is also bounded below on $\cup B_{n}$, then $y^{\ast}=x^{\ast}$;

\smallskip(e) for every unbounded nested sequence $\{B_{n}\}$ of balls such
that $x^{\ast}$ is bounded below on $\cup B_{n}$, $\cup B_{n}$ is an affine
half-space determined by $x^{\ast}$.

\medskip

\emph{Proof} (of \cite[Th.\ 6]{BaHuLiTr00} and \cite[Prop.\ 2.5]{BaHuLi04}
with $\lim_{m}\liminf_{n}$ instead of $\lim_{m}\lim_{n}$). We keep some parts
of the original proofs of \cite[Th.\ 6]{BaHuLiTr00} and \cite[Prop.\ 2.5]%
{BaHuLi04}.

Observe first that for $z\in X,$ $z^{\ast}\in X^{\ast}$ and $r>0$ one has
\begin{equation}
\inf z^{\ast}\left(  B(z,r\right)  )=\inf z^{\ast}(z+r\cdot\operatorname*{int}%
B_{X})=z^{\ast}(z)+r\inf z^{\ast}(B_{X})=z^{\ast}(z)-r. \label{rz0}%
\end{equation}

(a)~$\Rightarrow$ (b). Let $\{x_{n}^{\ast}\}\subseteq B_{X^{\ast}}$ and
$\{x_{m}\}\subseteq B_{X}$ be such that (\ref{rz1b}) is verified. As
$\{x_{n}^{\ast}\}\subseteq B_{X^{\ast}},$ $\{x_{n}^{\ast}\}$ has w*-cluster
points in $B_{X^{\ast}}$. Let $y^{\ast}$ be a w*-cluster point of
$\{x_{n}^{\ast}\}$. {Hence there exists a subnet $\{x_{n_{i}}^{\ast}\}$ {of}
$\{x_{n}^{\ast}\}$ {such that }$x_{n_{i}}^{\ast}\rightarrow^{w^{\ast}}y^{\ast
}$ $(\in B_{X^{\ast}})$, and so $\left\Vert y^{\ast}\right\Vert \leq1$ {and }%
\begin{align*}
\mu_{m}:=  &  \liminf_{n}\left(  \frac{x_{n}^{\ast}+x^{\ast}}{2}\right)
(x_{m})\leq\lim_{i}\left(  \frac{x_{n_{i}}^{\ast}+x^{\ast}}{2}\right)
(x_{m})=\left(  \frac{y^{\ast}+x^{\ast}}{2}\right)  (x_{m})\\
\leq &  \tfrac{1}{2}\left\Vert x_{m}\right\Vert \cdot\left\Vert y^{\ast
}+x^{\ast}\right\Vert \leq\tfrac{1}{2}\left\Vert y^{\ast}+x^{\ast}\right\Vert
\leq1\ \ \forall m\geq1.
\end{align*}
{By (\ref{rz1b}) one has }$\lim_{m}\mu_{m}=1$,{ whence }$\tfrac{1}%
{2}\left\Vert y^{\ast}+x^{\ast}\right\Vert =1${. Using (a) one gets }$y^{\ast
}=x^{\ast}${. This shows that the sequence }$\{x_{n}^{\ast}\}$ {has a unique
w*-cluster point }$x^{\ast}${, that is, w*-}$\lim y_{n}^{\ast}=x^{\ast}$.
Hence $x^{\ast}$ {is a w*-ALUR point of }$B_{X^{\ast}}$.}

\smallskip(b) $\Rightarrow$ (b$^{\prime}$). Suppose $\{y_{n}^{\ast}\}\subseteq
B_{X^{\ast}}$ is asymptotically normed by $B_{X}$ $(\subseteq B_{X^{\ast\ast}%
})$, where $y_{n}^{\ast}:=(x_{n}^{\ast}+x^{\ast})/2$ for $n\geq1$. By
definition, {for any $m\geq1$, there exists $x_{m}\in B_{X}${ and $N_{m}%
\in\mathbb{N}$ such that $(1\geq)$ $y_{n}^{\ast}(x_{m})>1-1/m$ for all $n\geq
N_{m}$. It follows that }$(1\geq)$ $\liminf_{n}y_{n}^{\ast}(x_{m})\geq1-1/m$
for $m\geq1$, and so{ $\lim_{m}\liminf_{n}y_{n}^{\ast}(x_{m})=1$.} Then
$y_{n}^{\ast}\rightarrow^{w^{\ast}}x^{\ast}$ by (b), whence $x_{n}^{\ast
}\rightarrow^{w^{\ast}}x^{\ast}$. }

{\smallskip}

(b$^{\prime}$) $\Rightarrow$ (c). Let $\{B_{n}:=B(x_{n},r_{n})\}$ be an
unbounded nested sequence of balls such that $x^{\ast}$ is bounded below on
$\cup B_{n},$ that is, $\{\inf x^{\ast}(B_{n})\}$ is bounded below; hence
{$0<r_{n}\uparrow\infty$}. Let $\{y_{n}^{\ast}\}\subseteq S_{X^{\ast}}$ be
such that $\{\inf y_{n}^{\ast}(B_{n})\}$ is bounded below, too. Take
$c\in\mathbb{R}$ a common lower bound. We may assume, without loss of
generality, that $0\in B_{1}$ $(${$\subseteq B_{n}$ for $n\geq1)$, and so
$\left\Vert 0-x_{n}\right\Vert =\left\Vert x_{n}\right\Vert <r_{n}$ for
$n\geq1$}. Set $y_{n}:=x_{n}/r_{n}$; it follows that $\left\Vert
y_{n}\right\Vert <1$. {Hence, by (\ref{rz0}), one has }%
\begin{equation}
c\leq\inf y_{n}^{\ast}(B_{n})=y_{n}^{\ast}(x_{n})-r_{n}\leq0\ \ \text{and
\ }c\leq\inf x^{\ast}(B_{n})=x^{\ast}(x_{n})-r_{n}\leq0\text{\ }\ \forall
n\geq1, \label{rz3}%
\end{equation}
{both inequalities \textquotedblleft$...\leq0$" being true because $0\in
B_{n}$. For $n\geq m$ one has $B_{n}\supseteq B_{m}$ (because $\{B_{n}\}$ is
increasing), whence $y_{n}^{\ast}(B_{n})\supseteq y_{n}^{\ast}(B_{m})$, and
so, by (\ref{rz3}) and (\ref{rz0}), one has
\[
c\leq y_{n}^{\ast}(x_{n})-r_{n}=\inf y_{n}^{\ast}(B_{n})\leq\inf y_{n}^{\ast
}(B_{m})=y_{n}^{\ast}(x_{m})-r_{m}\ \ \forall n,m\in\mathbb{N}^{\ast},\ n\geq
m.
\]
Hence $x^{\ast}(y_{m})\geq1+c/r_{m}$ for $m\geq1$ and $y_{n}^{\ast}(y_{m}%
)\geq1+c/r_{m}$ for $n\geq m\geq1$. It follows that%
\[
\forall m\geq1,\ \forall n\geq m\ :\left(  \frac{x^{\ast}+y_{n}^{\ast}}%
{2}\right)  (y_{m})\geq1+c/r_{m}.
\]
For }$\varepsilon>0$, there exists $m\geq1$ such $c/r_{m}>-\varepsilon$, and
so for $n\geq m$ one has that $\left(  \frac{x^{\ast}+y_{n}^{\ast}}{2}\right)
(y_{m})>1-\varepsilon$, proving so that $(B_{X^{\ast}}\supseteq)$ $\{(x^{\ast
}+y_{n}^{\ast})/2\}$ is asymptotically normed by $B_{X}$. By our hypothesis,
w*-$\lim\tfrac{1}{2}\left(  x^{\ast}+y_{n}^{\ast}\right)  =x^{\ast}$, and so
w*-$\lim y_{n}^{\ast}=x^{\ast}$.

\smallskip(c)~$\Rightarrow$ (d). Apply (c) to the constant sequence
$\{y_{n}^{\ast}=x^{\ast}\}$.

\smallskip(d)~$\Rightarrow$ (e). {Because $B_{n}$ is convex and open for
$n\geq1$ and $\{B_{n}\}$ is nested (hence increasing), $B:=$}$\cup B_{n}${ is
convex and open. Then, setting }$\alpha:=\inf x^{\ast}(B)$ $(\in\mathbb{R})$,
one has $B\subseteq\{x:x^{\ast}(x)>\alpha\}=:H$. If {$H\neq B$ $(\subseteq
H)$, there exists $z\in H\setminus B$; hence }$x^{\ast}(z)>\alpha$.{ By a
separation theorem there exists $y^{\ast}\in S_{X^{\ast}}$ such that $y^{\ast
}(z)\leq\inf y^{\ast}(B)$. Then }$y^{\ast}=x^{\ast}$ by our hypothesis, and so
$\alpha<x^{\ast}(z)=y^{\ast}(z)\leq\inf y^{\ast}(B)=\alpha$. Therefore, $\cup
B_{n}$ $(=H)$ is an open affine half-space determined by $x^{\ast}$.

\smallskip(e)~$\Rightarrow$ (a). Suppose there exists $y^{\ast}\in S_{X^{\ast
}}\setminus\{x^{\ast}\}$ such that $\tfrac{1}{2}(x^{\ast}+y^{\ast})\in
S_{X^{\ast}}$. Let $\{\delta_{n}\}\subseteq(0,1)$ such that $\sum
_{n=1}^{\infty}\delta_{n}<1$ and $\{x_{n}\}\subseteq B_{X}$ such that
$(x^{\ast}+y^{\ast})(x_{n})>2-\delta_{n}$ for $n\geq1$. Then $(0+0\leq)$
$[1-x^{\ast}(x_{n})]+[1-y^{\ast}(x_{n})]=2-(x^{\ast}+y^{\ast})(x_{n}%
)<\delta_{n}$, and so $x^{\ast}(x_{n})-1>-\delta_{n}$ and $y^{\ast}%
(x_{n})-1>-\delta_{n}$ for $n\geq1$.

Clearly, $r_{n}:=n+\sum_{i=1}^{n}\delta_{i}$ $(\uparrow\infty)$. Set $u_{n}:=
{\textstyle\sum\nolimits_{i=1}^{n}} x_{i}$ and $B_{n}:=B(u_{n},r_{n})$. Take
$x\in B_{n}$; then
\[
\left\Vert x-u_{n+1}\right\Vert =\left\Vert x-(u_{n}+x_{n+1})\right\Vert
\leq\left\Vert x-u_{n}\right\Vert +\left\Vert x_{n+1}\right\Vert \leq
n+{\textstyle\sum\nolimits_{i=1}^{n}} \delta_{i}+1<r_{n+1},
\]
and so $x\in B_{n+1}$. Hence $B_{n}\subseteq B_{n+1}$, and so $\{B_{n}\}$ is
an unbounded nested sequence of balls.

For any $n\geq1,$%
\[
\alpha:=\inf x^{\ast}(B_{n})=x^{\ast}({\textstyle\sum\nolimits_{i=1}^{n}}%
x_{i})-n-{\textstyle\sum\nolimits_{i=1}^{n}}\delta_{i}={\textstyle\sum
\nolimits_{i=1}^{n}}[x^{\ast}(x_{i})-1-\delta_{i}]>-{\textstyle\sum
\nolimits_{i=1}^{n}}2\delta_{i}>-2.
\]
Similarly, $\beta:=\inf y^{\ast}(B_{n})>-2$. Hence both $x^{\ast}$ and
$y^{\ast}$ are bounded below on $B=\cup B_{n}$. Because $B$ is an open and
affine half-space determined by $x^{\ast}$ and $y^{\ast}$, one obtains that
\[
B=\{x\mid x^{\ast}(x)>\alpha\}=\{x\mid y^{\ast}(x)>\beta\}.
\]
Fixing some $x_{0}\in B$, $u\in\lbrack x^{\ast}\geq0]$ and $t\geq0$
one has $x_{0}+tu\in B,$ and so
$\beta<y^{\ast}(x_{0})+ty^{\ast}(u)$. As $t\geq0$ is arbitrary one
has that $y^{\ast}(u)\geq0$, and so $x^{\ast}(u)\geq0\Rightarrow
y^{\ast}(u)\geq0$. It follows that $\ker x^{\ast}\subset\ker
y^{\ast}$, and so $y^{\ast}=\gamma x^{\ast}$ for some
$\gamma\in\mathbb{R}$.  Since $x^*, y^* \in S_{X^*}$ we have $y^*=
x^*$ which is a contradiction. \hfill$\square$

\medskip Notice that we did not check other proofs from the cited works.


\end{document}